\documentclass[12pt,twoside]{gen-j-l}
\usepackage{amsthm}
\usepackage{amsmath}
\usepackage{amstext}
\usepackage{amsfonts}
\usepackage{latexsym}
\usepackage{amscd}
\newcommand{\la}{\ensuremath{\longrightarrow}}
\newcommand{\Y}{\ensuremath{\mathrm{Y}}}

\newcommand{\hilb}{\ensuremath{\mathrm{Hilb}}}
\newcommand{\p}{\ensuremath{\mathbb{P}}}
\newcommand{\pone}{\ensuremath{\mathbb{P}^{1}}}
\newcommand{\s}{\ensuremath{\mathrm{S}}}

\newcommand{\x}{\ensuremath{\mathrm{X}}}
\newcommand{\Z}{\ensuremath{\mathrm{Z}}}
\newcommand{\C}{\ensuremath{\mathrm{C}}}

\newcommand{\I}{\ensuremath{\mathcal{I}}}

\newcommand{\sheaf}{\ensuremath{\mathcal{O}}}

\newcommand{\E}{\ensuremath{\mathcal{E}}}
\newcommand{\A}{\ensuremath{\mathbf{A}}}
\newtheorem{theorem}{Theorem}[section]
\newtheorem{lema}[theorem]{Lemma}
\newtheorem{definition}[theorem]{Definition}
\newtheorem{proposition}[theorem]{Proposition}
\newtheorem{corollary}[theorem]{Corollary}

\begin{document}
\title[Extensions of \pone\ and their Hilbert schemes]{Infinitesimal extensions of \pone\ and their Hilbert schemes}
\author{Nikolaos Tziolas}
\address{Deparment of Mathematics, University of Utah, Salt Lake City, Utah, 84112.}
\address{\textit{Current address:}\textrm{ Mathematics Institute, University of Warwick, Coventry, CV4 7AL, United Kingdom}}
\email{tziolas@maths.warwick.ac.uk}
\thanks{I would like to thank my advisor J\'{a}nos Koll\'{a}r without whose valuable advice and 
continuous support this work wouldn't have been possible.}





\begin{abstract}
In order to calculate the multiplicity of an isolated rational curve $C$ on a local complete intersection variety $X$, i.e. the 
length of the local ring of the Hilbert Scheme of $X$ at $[C]$, 
it is important to study infinitesimal neighborhoods of the curve in $X$. This is equivalent to infinitesimal 
extensions of \pone\ by locally free sheaves. In this paper we study infinitesimal extensions of \pone, determine 
their structure and obtain upper and lower bounds for the length of the local rings of their Hilbert schemes at [\pone].
\end{abstract}

\maketitle

\setcounter{section}{-1}
\section{Introduction}

The problem of counting curves of a certain ``type'' on an algebraic variety is a very old and difficult one. 
It is classically known that there are exactly 27 distinct lines on a smooth cubic surface in $\p^3$ and 2875 on a general 
quintic 3-fold in $\p^4$. Recent advances in string theory and mirror symmetry revived the problem of counting rational curves 
on Calabi-Yau 3-folds and have thrown new light in it. In particular it is of interest to know the contribution of an 
isolated rational curve $C$ to the total number of curves in a Calabi-Yau 3-fold $X$. In~\cite{Tziolas99} this 
contribution, the multiplicity of the curve, was defined to be the length of the local ring of the Hilbert scheme of $X$ at $[C]$, the point 
corresponding to $C$, and an explicit algorithm and formula was given under certain semipositivity conditions 
on $\I_{C,X} / \I_{C,X}^{(2)}$.

The problem is also of interest from the point of view of birational geometry. If $X \la Y$ is a birational map of 3-folds, with 
$X$ smooth, contracting a single curve $C$, then it is known that $C \cong \pone$ and $\I / \I^2 $ is isomorphic to one of 
$\sheaf_{\pone}(1) \oplus  \sheaf_{\pone}(1)$, $\sheaf_{\pone} \oplus \sheaf_{\pone}(2)$ or $\sheaf_{\pone}(-1) \oplus \sheaf_{\pone}(3)$. 
In particular this is the case of a flopping contraction that appears in the Minimal Model Program or if $X$ is Calabi-Yau. 
One would like to study such curves even when $X$ has terminal singularities. A $(1,1)$ curve is rigid and it contracts. 
The $(0,2)$ case was studied by Miles Reid~\cite{Reid83} and it either contracts or it moves in a positive dimensional family. 
The method used to prove this, as in~\cite{Tziolas99} was the study of infinitesimal neighborhoods of the curve in $X$ which 
is the same as infinitesimal extensions of \pone\ by $\sheaf_{\pone}$. 

The case of $(-1,3)$ curves is still open. There are examples by Laufer and Jimenez~\cite{Jim92} of $(-1,3)$ curves that 
neither contract or move. To study this case, or calculate the multiplicity of a curve in a Calabi-Yau without any semipositivity conditions, 
one can use the same technique as before, i.e. study infinitesimal neighborhoods of the curve or equivalently infinitesimal 
extensions of \pone with locally free sheaves.

In this paper i want to remove the semipositivity conditions that appear in~\cite{Tziolas99} and following the ideas exposed there 
study infinitesimal extensions of \pone\ with seminegative locally free sheaves and in particular their Hilbert schemes. 
The situation here is a lot more subtle than the case of extensions by $\sheaf_{\pone}$. Example 1.1 shows that it is no longer true 
that the number of extensions determine the Hilbert scheme~\cite{Tziolas99}. It even depends on the characterstic of the base field. 

In section 1, we use the theory of Hochschild extensions of commutative algebras~\cite{Weibel94} to describe the structure of infinitesimal 
extensions of \pone\ by locally free sheaves. 

In section 4, there is an example, suggested to me by J\'{a}nos Koll\'{a}r, 
of a rigid smooth rational curve $C$ in a surface $S$ with $\I^{(n)}/\I^{(n+1)} \cong \sheaf_{\pone}(-1)$ $
\forall n \geq 1$. This shows that a result similar to~\cite[theorem 3.1]{Tziolas99} giving an explicit formula 
to calculate the multiplicity of a curve is unlikely to exist. So we settle for finding 
upper and lower bounds for the length of the local ring of the Hilbert scheme at $[\pone]$ of a scheme \Z\ obtained by extending \pone. 
Moreover, this example shows that a lower bound other than a constant cannot exist.

In section 2, theorem 2.1 shows that any extension can be deformed to the trivial one, and theorem 2.2 calculates Hilbert schemes of 
trivial extensions, and hence establishes an upper bound. We also show by examples that the length can differ from the length 
of the corresponding trivial extension.

Finally in section 3, proposition 3.2 establishes a lower bound.

\section{Structure of the extensions}
Let \Z\ be a scheme obtained from \pone\ by a sequence of infinitesimal extensions by $\sheaf_{\pone}$. 
Then by~\cite[Proposition 4.1]{Tziolas99}, 
$Z \cong \pone \times Spec D $, and $D=H^{0}(Z,\sheaf_{Z})$. Hence $\sheaf_{\hilb(Z),[\pone]}
=D$. In particular $length_{[\pone]} \hilb(Z)$ depends 
only on the number of extensions. This is no longer true in the case of extensions by 
$\sheaf_{\pone}(-1)$ as seen by the next example.

\textbf{Example 1.1}: In this example I will exhibit two schemes \x\ and $Y$, both obtained as infinitesimal extensions 
of \pone\ by two $\sheaf_{\pone}(-1)$, 
such that $length_{[\pone]} \hilb(X)=6$ and $length_{[\pone]} \hilb(Y)=5$.

(a) Let $X=2L$ be a double line in $\p^{3}$, i.e if $\I$ is the ideal of a line, then \x\ is defined 
by $\I^{2}$. Let $x_{0},\; x_{1}, \; x_{2},
\;x_{3}$ be the coordinates in $\p^{3}$ and $\I=(x_{0},x_{1})$. A line $L$ in  $\p^{3}$ has parametric equations \[
\begin{array}{c}
x_{0}=ax_{2}+bx_{3} \\
x_{1}=cx_{2}+dx_{3}
\end{array} \]
Clearly $\I/\I^{2} \cong \sheaf_{\pone}(-1) \oplus \sheaf_{\pone}(-1)$, and hence 
\x\ is an extension of \pone\ by two $\sheaf_{\pone}(-1)$. The condition for a line $L$ to be in \x\ is \[
\begin{array}{ccc}
(ax_{2}+bx_{3})^{2}=0 & (cx_{2}+dx_{3})^{2}=0 & (ax_{2}+bx_{3})(cx_{2}+dx_{3})=0.
\end{array} \]
From these equations we immediately see that \[
\sheaf_{\hilb(X),[L]}= \frac{\textstyle k[a,b,c,d]}{\textstyle (a^{2}, b^{2}, c^{2}, d^{2}, 2ab, 2cd, ac,bd,ad+bc)} \]
which has length 6 if $ch(k)\neq 2$ and 8 if $ch(k)=2$. This is by itself surprising since one would expect that the first order 
neighborhood of a line contains exactly the first order deformations, i.e if $G=G(2,4)$ is the grassmanian 
of lines in $\p^{3}$ and $
[L]$ the point corresponding to the line, then one may expect that $ length \sheaf_{\hilb(2L),[L]}= length \sheaf_{G,[L]} / m_{[L]}^{2}=5 $.

(b) Let $S \subset \p^{3}$ be the cubic surface with a $D_{5}$ singularity given by the equation 
$f=x_{3}x_{0}^{2}+x_{0}x_{2}^{2}+x_{2}x_{1}^{2}$, and let $L=(x_{0},x_{1})$, a line through the singular 
point. I will show that 
$\I / \I^{(2)} \cong \sheaf_{L}(-1)$, $\I^{(2)} / \I^{(3)} \cong \sheaf_{L}(-1)$, and 
 $length \sheaf_{\hilb(SpecY),[L]}=5$, where $Y=Spec(\sheaf_{S} / \I^{(3)})$ appears as an 
infinitesimal extension of \pone\ by two $\sheaf_{L}(-1)$. i.e 
\begin{gather*}
0\la \sheaf_{L}(-1)=\I / \I^{(2)} \la \sheaf_{S} / \I^{(2)} \la \sheaf_{S} / \I = \sheaf_{L} \la 0 \\
0\la \sheaf_{L}(-1)=\I^{(2)} / \I^{(3)} \la \sheaf_{S} / \I^{(3)} \la \sheaf_{S} / \I^{(2)} \la 0.
\end{gather*}
Then \[
\I^{2}= \frac{\textstyle (x_{0}^{2},x_{1}^{2},x_{0}x_{1},f)}{\textstyle(f)}=
\frac{\textstyle (x_{0}^{2},x_{1}^{2},x_{0}x_{1},x_{0}x_{2}^{2})}{\textstyle (f)}. \]
Since $x_{0}x_{2}^{2} \in \I^{2}$ and $x_{2}^{2} \notin \I$, then $x_{0} \in \I^{(2)}$. Hence \[
\I^{(2)}=(x_{0},x_{1}^{2})/(f). \]
Now the map \[
\I / \I^{(2)}=(x_{0},x_{1}) / (x_{0},x_{1}^{2}) \ni x_{0}h+x_{1}g \longmapsto g(0,0,x_{2},x_{3}) \in k[x_{2}
,x_{3}] \]
is an isomorphism that shows that \[
 \I / \I^{(2)} \cong \sheaf_{L}(-1). \]
Moreover, $\I \I^{(2)} = (x_{0}^{2},x_{0}x_{1},x_{1}^{3},x_{0}x_{2}^{2}+x_{2}x_{1}^{2})/(f)$. Since 
$x_{2}(x_{0}x_{2}+x_{1}^{2}) \in \I / \I^{(2)}$, then $x_{0}x_{2}+x_{1}^{2} \in \I^{(3)}$. Hence \[
\I^{(3)}=(x_{0}^{2},x_{1}^{3},x_{0}x_{1},x_{0}x_{2}+x_{1}^{2})/(f). \]
The map \[
\I^{(2)} / \I^{(3)} \ni x_{0}h + x_{1}^{2}g \longmapsto h(0,0,x_{2},
x_{3})-x_{2}g(0,0,x_{2},x_{3}) \in k[x_{2},x_{3}] \]
is an isomorphism which shows that \[ 
\I^{(2)} / \I^{(3)} \cong \sheaf_{L}(-1). \]
Doing a calculation as in part $(a)$, one can easily see that \[
\sheaf_{\hilb(Y),[L]} = \frac{\textstyle k[a,b,c,d]}{\textstyle (a^{2},b^{2},c^{3},d^{2}, 2ab,3c^{2}d,3d^{2}ac,bc,ad+bc,a+c^{2},b+2cd)} \]
which has length 5 if $ch(k) \neq 2,3$. \qed

In section 4 there is an even more startling example of how strange extensions by $\sheaf_{\pone}(-1)$ can 
really be.

So a result similar to~\cite[Theorem 3.1]{Tziolas99} is rather unlikely to exist, at least not 
without further assumptions. What we will do next is to obtain upper and lower bounds for the $length \sheaf_{\hilb(X),[\pone]}$, 
\x\  obtained by 
extending \pone\ by $\sheaf_{\pone}(-d)$'s. 

The theory of Hochschild extensions of commutative algebras~\cite{Weibel94} can be used to understand 
the structure of such extensions. Arguing as in the ring case we get the following
\begin{lema}\label{ringstructure}
Let $Y \stackrel{h}{\la} Z$ be a morphism of schemes having the same underlying topological space. 
 Let \x\ be an infinitesimal extension of $Y$ by a coherent 
sheaf $\mathcal{F}$, such that $h$ factors through a morphism $X \la Z$, i.e $\sheaf_{X}$ is an $\sheaf_{Z}$-algebra, and 
the extension is split as an extension of $\sheaf_{Z}$-modules. Then of course $\sheaf_{X} \cong \sheaf_{Y} \oplus \mathcal{F}$ 
as an $\sheaf_{Z}$-modules, and the ring structure is defined by a map $f \in \mathrm{Hom}_{Z}(\sheaf_{Y} \otimes_{\sheaf_{Z}} 
\sheaf_{Y} , \mathcal{F})$, such that
\begin{enumerate} 
\item If $U \subset Y$ is open affine, then \[
r_{0}f(r_{1},r_{2})-f(r_{0}r_{1},r_{2})+f(r_{0},r_{1}r_{2})-f(r_{0},r_{1})r_{2}=0 \]
for all $r_{i} \in \sheaf_{Y}(U)$.
\item $f(r_{0},r_{1})=f(r_{1},r_{0})$, $ \forall r_{i} \in \sheaf_{Y}(U)$, i.e $f$ is symmetric.
\item The multiplication is defined by \[
(x,m)(x^{\prime},m^{\prime})=:(xx^{\prime},xm^{\prime}+x^{\prime}m +f(x,x^{\prime})) \]
for $x \in \sheaf_{Y}(U) $ and $m \in \mathcal{F}(U)$. It is associative by 1, commutative by 2 and the unit is $
(1,-f(1,1))$.
\end{enumerate}
\end{lema}
Hochschild cohomology~\cite{Weibel94} can be naturally defined for ringed spaces and can be used to classify extensions 
as in the previous lemma. With assumptions as in the lemma, consider the complex
\[
\mathrm{Hom}_{Z}(\sheaf_{Y} , \mathcal{F}) \stackrel{\delta^{0}}{\la} 
\mathrm{Hom}_{Z}(\sheaf_{Y}^{\otimes 2}, \mathcal{F}) 
\stackrel{\delta^{1}}{\la} \mathrm{Hom}_{Z}(\sheaf_{Y} ^{\otimes 3}, \mathcal{F}) \]
where for any open affine $U \subset Y$, $\delta^{0}$, $\delta^{1}$ are defined by the rule
\[
\delta^{0}(h)(r_{0} , r_{1}) = r_{0} h(r_{1})- h(r_{0}  r_{1}) + r_{1} h(r_{0} ) 
\]
for $h \in \mathrm{Hom}_{Z}(\sheaf_{Y} , \mathcal{F})$, and  \[
\delta^{1}(f)(r_{0} , r_{1} , r_{2})= 
r_{0} f(r_{1} , r_{2}) - f(r_{0} r_{1} , r_{2}) + f(r_{0} , r_{1} r_{2}) - f(r_{0} , r_{1}) r_{2} \]
for $f \in \mathrm{Hom}_{Z}(\sheaf_{Y} \otimes \sheaf_{Y} , 
\mathcal{F}) $. 

Now as in the ring case it is easy to see that two maps $f_{1}, \; f_{2}$ define the 
same ring structure if and only if $\exists \; h \in 
\mathrm{Hom}_{Z}(\sheaf_{Y} , \mathcal{F}) $ such that $
f_{1} - f_{2} = \delta^{0}(h) $. 
Hence if $H^{2}(\sheaf_{Z} , \mathcal{F})=Ker \delta^{1}/ Im \delta^{0}$, then: 
\begin{lema}\label{hochschild}
The space of isomorphism classes of infinitesimal extensions of $Y$ by $
\mathcal{F} $, that are also extensions of $\sheaf_{Z}$-algebras and split as extensions of $\sheaf_{Z}$-modules, 
is in one to one correspondence with the subspace 
$H^{2}_{S}( \sheaf_{Y} , \mathcal{F})$ of 
$H^{2}(\sheaf_{Y} , \mathcal{F}) $  consisting of classes of symmetric maps $f$.
\end{lema}

In particular, extensions of \pone\ by $\sheaf_{\pone}(d_{i})$, $ -3 \leq d_{i} $ are by~\cite[Proposition 4.1]{Tziolas99} 
extensions of the above type and hence lemma~\ref{hochschild} applies.

\section{Upper bound}

The next theorem compares any extension with the trivial one and is fundamental for obtaining upper bounds for the 
Hilbert schemes.

\begin{theorem}\label{deformtotrivial}
Let $Z$ be a scheme, $Y \stackrel{h}{\la} Z$ a scheme over $Z$ such that $Y_{top}=Z_{top}$, and 
there is a section $\sigma$ of $h$. (In other words $\sheaf_{Y}$ is an $\sheaf_{Z}$-algebra and $Z$ is also a 
closed subscheme of $Y$. Moreover, $Y$ is necessarily obtained from $Z$ by a sequence of infinitesimal extensions.)
Let \x\ be an infinitesimal extension of $Y$ by a coherent sheaf $\mathcal{F}$, i.e  \[
0 \la \mathcal{F} \la \sheaf_{X} \la \sheaf_{Y} \la 0 \]
such that this is an  extension of $\sheaf_{Z}$-algebras and split as extensions of $\sheaf_{Z}$-modules. Then 
\begin{enumerate}
\item $X$ can be deformed to the trivial extension $W$ of $Y$ by $\mathcal{F}$.
\item \[
length _{[Z]}\hilb(X) \leq length_{[Z]} \hilb(W) \]
\end{enumerate}
\end{theorem}

\begin{proof}
The assumptions of the theorem imply that lemma~\ref{ringstructure} applies and hence $\sheaf_{X} , \; \sheaf_{W} 
=\sheaf_{Y} \oplus \mathcal{F}$ as  $
\sheaf_{Z}$-modules, and the ring structure is defined by a map $f \in \mathrm{Hom}_{Z}(\sheaf_{Y} \otimes_{\sheaf_{Z}} 
\sheaf_{Y} , \mathcal{F})$ for $X$, satisfying the conditions of lemma~\ref{ringstructure}, and the zero map for $W$.
The idea of the proof is to try to deform $f$ to 0.

Consider the projection $Y  \times \A^{1} \stackrel{\pi}{\la} Y$. 
Construct an infinitesimal extension \s, of $Y \times \A^{1}$ by $\pi^{\ast} \mathcal{F}$  \[
0 \la \pi^{\ast} \mathcal{F} \la \sheaf_{S} \la \sheaf_{Y \times \A^{1}} \la 0 \]
by letting $\sheaf_{S}= \sheaf_{Y \times \A^{1}} \oplus \pi^{\ast} \mathcal{F} $ as a group, and 
define it's ring structure by an $\sheaf_{Z \times \A^{1}}$-linear map \[
\sheaf_{Y \times \A^{1}} \otimes_{\sheaf_{Z \times \A^{1}}} \sheaf_{Y \times \A^{1}} \stackrel{\Phi}{\la} \pi^{\ast} \mathcal{F} \]
satisfying the conditions of lemma~\ref{ringstructure}. 
Construct this map locally. Let $U=SpecA \subset Y$ be affine open, and $
Z \cap U = SpecB$. Assume $\mathcal{F} \mid_{U} \cong \widetilde{M_{U}}$. 
Then by the assumptions on $Y$, $B \subset A$ and there is a surjection $A \la B$, 
such that the composition $B \la A \la B $ is the identity. 
Now  $U \times \A^{1} = SpecA[t] $, and $\pi^{\ast} \mathcal{F} \mid_{U \times \A^{1}} =  \widetilde{M_{U} \otimes_{A} A[t]}$. 
Working over $U$, we need to define a map \[
\Phi_{u}:A[t] \otimes_{B[t]} A[t] \la M_{U} \otimes_{B} B[t]. \]
Let $\phi(t)=\sum_{i}a_{i}t^{i} , \; h(t)=\sum_{i}a_{i}^{\prime}t^{i} \in A[t]$. Define $\Phi_{u}$ 
by\[
\Phi_{u}(\phi(t),h(t))=:t \sum_{i,j}f(a_{i},a_{j}^{\prime}) \otimes t^{i+j}. \]
Since $f$ satisfies the conditions of lemma~\ref{ringstructure}, then so does $\Phi_{u}$ and hence it defines a ring 
structure. Since $f$ is defined globally, the maps $\Phi_{u}$ glue to a global map $\Phi$ with the right 
properties. Now it is not difficult to check that $S \la \A^{1} $ is a deformation of $X$ to $X_{a}$, where $X_{a}$ is 
the infinitesimal extension of $Y$ by $\mathcal{F}$ corresponding to the map $af$. In particular for $a=0$ we get the trivial extension. 
Moreover since $X_{a} \cong X_{b}$ (as schemes not as extensions) $\forall a, b \neq 0$, part \textit{2} follows immediately.

\end{proof}
The next corollary is important for the applications. 
\begin{corollary}
With \Y\ and $Z$ as in the previous theorem, let $\x_{n}$ be obtained from \Y\ by succesive infinitesimal extensions by $\mathcal{F}_{i}$, $
i=1, \ldots , n$, such that they are extensions of $\sheaf_{Z}$-algebras and split as extensions of $\sheaf_{Z}$-modules. Then \[
length_{[Z]} \hilb(\x_{n}) \leq length_{[Z]} \hilb(W) \]
where $W$ is the trivial infinitesimal extension of $Z$ by $\oplus_{i=1}^{n} \mathcal{F}_{i}$.
\end{corollary}

\begin{proof}
This is proved by applying theorem~\ref{deformtotrivial}  many times. Start with $\x_{1}$. This is an extension 
of \Y\ by $\mathcal{F}_{1}$. Let $\x_{1}^{\prime}$ be the trivial extension of $Y$ by 
$\mathcal{F}_{1}$. Now  let $\x_{2}^{\prime}$ be the trivial extension of $\x_{1}$ by $\mathcal{F}_{2}$. Then it is easy to see that 
it appears as a (possibly nontrivial) extension \[
0 \la \mathcal{F}_{2} \la \sheaf_{\x_{2}^{\prime}} \la \sheaf_{\x_{1}^{\prime}} \la 0. \]
Now if $\x_{2}^{\prime \prime}$ is the trivial extension of $\x_{1}^{\prime}$ by $\mathcal{F}_{3}$, we 
get from the theorem that \[
length_{[Z]} \hilb(\x_{2}^{\prime}) \leq length_{[Z]} \hilb(\x_{2}^{\prime \prime}) \] 
and hence \[ 
length_{[Z]} \hilb(\x_{2}) \leq length_{[Z]} \hilb(\x_{2}^{\prime \prime}). \] 
Continuing this way we get the corollary.
\end{proof}
So in order to get upper bounds for the length of the Hilbert schemes we only need to 
study split infinitesimal extensions  of \pone. This is done in the following theorem.

\begin{theorem}
Let $\E ^{\prime} = Ample \oplus \E$, with $\E = \oplus_{i=1}^{s} \sheaf_{\pone}(-d_{i})$, $d_{i} \geq 0$, and 
let $X$ be a trivial infinitesimal extension of \pone\ by $\E ^{\prime}$. Let $D = \sheaf_{\hilb(X),[\pone]}$, $m$ the maximal ideal. Then
\begin{enumerate} 
\item If the characteristic of the base field is different than 2, then:
\begin{enumerate} 
\item If $\mathcal{E}=\sheaf_{\pone}(-d) \;\; , d \geq 0$ then we get the total length
\[
length_{[\pone]} \hilb(X) = \sum_{i \geq -2} \left( \begin{array}{c}
d-i \\
2+i
\end{array} \right). \]
\item  For any \E, \[
length( D/ m^{3}) = 1 + (r(\mathcal{E}) - \deg \mathcal{E}) + \left( \begin{array}{c}
-\deg \E \\
2
\end{array} \right).  \]
If $0 \leq d_{i} \leq 3$ then $m^{3}=0$ and the formula gives the total length.
\item \[
length( D/ m^{4})=  length( D/ m^{3}) + 
 \left( \begin{array}{c}
-\deg \E - r(\E) \\
3
\end{array} \right) - \]
\[
- (r(\E)-1)(r(\E)^{2}-2r(\E)-9r(\E) \deg \E -18 \deg \E ). 
\]
If $0 \leq d_{i} \leq 5$, then $m^{4}=0$ and the formula gives the total length.
\end{enumerate}
\item If the characteristic of the base field is 2 and $\mathcal{E}=\sheaf_{\pone}(-d) \; d \geq 0$,
 then \[
\sheaf_{\hilb(X),[\pone]} = \frac{\textstyle k[z_{0},\ldots , z_{d}]}{\textstyle  (z_{0}^{2},\ldots ,z_{d}^{2})} \]
and hence \[
length_{[\pone]} \hilb(X) = 2^{d+1}. \]
\end{enumerate}
\end{theorem}
For ease in applications we get:
\begin{corollary}
With assumptions as in the previous theorem,
\begin{enumerate}
\item If $ \E = \oplus_{i=1}^{n} \sheaf_{\pone}(-1) $, then \[
\sheaf_{\hilb(X),[\pone]} = \frac{\textstyle k[x_{1}, y_{1},x_{2}, y_{2},\ldots ,x_{n}, y_{n}]}{\textstyle ( x
_{i}x_{j}, 
y_{i}y_{j}, x_{i}y_{j}+x_{j}y_{i}, \;\; /1 \leq i \leq j \leq n )} \]
and \[
length_{[\pone]} \hilb(X) = 2n+1+ \frac{\textstyle n(n-1)}{\textstyle 2}. \]
Hence if $L$ is a line in $\p ^{n+1}$, then $length_{[L]} \hilb(2L)=2n+1+n(n-1)/2$.
\item If $\E = \oplus_{i=1}^{n} \sheaf_{\pone}(-2) $, then \[
\begin{array}{c}
\sheaf_{\hilb(X),[\pone]}= \\
\frac{\textstyle k[y_{1}, z_{1}, w_{1}, \dots , y_{n}, z_{n}, w_{n}]}{\textstyle ( y_{i}y_{j}
, w_{i}w_{j},z_{i}w_{j}+z_{j}w_{i},y_{i}z_{j}+y_{j}z_{i}, y_{i}w_{j}+y_{j}w_{i}+z_{i}z_{j}, \;\; 
/1 \leq i \leq j \leq n )} 
\end{array} \]
and hence \[
length_{[\pone]} \hilb(X) =3n+1+n(2n-1). \]
\end{enumerate}
\end{corollary}
Note that part $1$ of the corollary shows that in the case of a line in $\p^{n+1}$ the first order neighborhood does 
not contain only the 2n+1 
first order deformations in $\p^{n+1}$ as one would expected, but it also contains $n(n-1)/2$ second order deformations which is 
surprising. The case of a conic in $\p^{2}$ is even more surprising. Part $2$ shows that the Hilbert scheme of the first order 
neighborhood has length 13=5+7+1, and hence the first order neighborhood not only contains the 5 first order deformations, but 
also 7 second and 1 third order deformation, which is unexpexted.
\begin{proof}[Proof of Theorem 2.3]
Since $X$ is a trivial extension, it is nothing but the first order neighborhood of the section of the bundle 
$\p (\sheaf_{\pone} \oplus \E ) \la \pone$ corresponding to the first projection $
\sheaf_{\pone} \oplus \E \la \sheaf_{\pone} $. Let $z_{i,j}$, $0 \leq i \leq d_{j}$, $1\leq j \leq s$ 
be the local coordinates of $\hilb(X)$ at $[\pone]$. They correspond to 
a section of the bundle given locally by \[
t \longmapsto (t, \sum_{i=0}^{d_{1}}z_{i,1}t^{i} , \ldots , \sum_{i=0}^{d_{s}}z_{i,s}t^{i} ). \] 
The condition that this
 stays in the first order neighborhood of the section corresponding to the first projection as above is that \[
(z_{0,m}+z_{1,m}t + \ldots z_{d_{m},m}t^{d_{m}})  
(z_{0,\lambda}+z_{1,\lambda}t + \ldots z_{d_{\lambda},\lambda}t^{d_{\lambda}}) =0. \]
Now it immediately follows that the equations of $\hilb(X)$ at $[\pone]$ are 
\begin{equation}
 \sum_{i+j=\nu} z_{i,m}z_{j,\lambda}=0. 
\end{equation}
I will only give the proof of the case that $\E = \sheaf_{\pone}(-d)$. This is the only case that we get a nice formula. 
The other results stated are proved similarly. 

Let $S_{d}=:k[z_{0}, \ldots , z_{d}]$, $I_{d}=:(f_{0},\ldots , f_{2d})$, where $f_{m}=\sum_{i+j=m}z_{i}z_{j}$, and $A_{d}=:S_{d}/I_{d}$. 
The key point of the proof is that \\
\textbf{Claim 1}. \emph{The elements $f_{m}$ form a Gr\"{o}bner basis for $I_{d}$ with respect to reverse lexicographic order}.

From now on we work with reverse lexicographic order. For any $f \in S_{d} $, let $in(f)$ be the initial part of $f$ with respect to the order. 
The following results from the theory of Gr\"{o}bner bases are needed.
\begin{definition}[\protect{~\cite[Definition 15.6]{Eis94}}] 
Let $S$ be a polynomial ring with monomial order $>$. If $f, g_{1}, \ldots , g_{t} \in S$ then there is an expression \[
\begin{array}{ccc}
f=\sum_{i}f_{i}g_{i} + f^{\prime} & with \; f^{\prime} \in S & f_{i} \in S 
\end{array} \]
where none of the monomials of $f^{\prime}$ is in $(in(g_{1}), \ldots , in(g_{t}))$ and \[
in(f) \geq in(f_{i}g_{i}) \]
for every i. Any such $f^{\prime}$ is called a remainder of $f$ with respect to $ g_{1}, \ldots , g_{t}$, 
and an expression $f=\sum_{i}f_{i}g_{i} + f^{\prime}$ satisfying the above conditions is called a standard expression 
for $f$ in terms of the $g_{i}$.
\end{definition}
\begin{theorem}[Buchberger's Criterion, \protect{\cite[Theorem 15.8]{Eis94}}]
Let\\ $S$ be a polynomial ring with monomial order $>$, and $I \subset S$ an ideal. Let $g_{1},\ldots,g_{t} \in I$, and $
m_{ij}=in(g_{i})/GCD(in(g_{i}),in(g_{j}))$. Choose a standard expression \[
m_{ji}g_{i}-m_{ij}g_{j} = \sum_{u}f_{u}^{(ij)}g_{u} + h_{ij} \]
for $m_{ji}g_{i}-m_{ij}g_{j} $ with respect to $g_{1}, \ldots, g_{t}$. 

Then the elements $g_{1}, \ldots, g_{t}$ form a Gr\"{o}bner basis for $I$ iff $h_{ij}=0$ for all $i$ and $j$.
\end{theorem}
\begin{corollary}[\protect{~\cite[exercise 15.20]{Eis94}}] 
In the case of the previous theorem, we only need to check Buchberger's criterion for the elements 
$g_{i}$, $g_{j}$ with $GCD(in(g_{i}),in(g_{j})) \neq 1$.
\end{corollary}
To prove claim 1, we will use Buchberger's criterion. Assume that $ch(k) \neq 2$. Then $I_{d}=(f_{2m+1},2f_{2n} \; / \; n,m)$. 
The introduction of 2 is only to reduce the number of calculations needed. It is easy to see that \[
in(f_{k})=\left\{ \begin{array}{cc}
z_{k/2}^{2} & \mbox{\textit{if k is even}} \\
z_{k-1/2}z_{k+1/2} & \mbox{\textit{if k is odd}}
\end{array} \right. \]
We only need to check Buchberger's criterion if $GCD(in(f_{i}),in(f_{j})) \neq 1$. It is easy to see that there are only two cases.

\textbf{Case 1}: $i=2m+1$ and $j=2n$. Then $in(f_{2m+1})=z_{m}z_{m+1}$ and $in(2f_{2n})=2z_{n}^{2}$. The only case to check is when $n=m$ or 
$n=m+1$. So assume that $n=m$. The other case is similar. Then $GCD(in(f_{2n+1}),in(2f_{2n}))=z_{n}$, and \[
m_{2n,2n+1}=2z_{n}f_{2n+1}-z_{n+1}f_{2n} \]
and hence since $in(m_{2n,2n+1})=z_{n}^{2}z_{n+1}$, this is a standard expression and Buchberger's criterion applies.

\textbf{Case 2}: $i=2n+1$ and $j=2m+1$. Then $in(f_{2n+1})=z_{n}z_{n+1}$ and $in(f_{2m+1})=z_{m}z_{m+1}$. The only case to check 
is if $n=m+1$ or $m=n+1$. Assume that $n=m+1$. Then $GCD(in(f_{2m+1}),in(2f_{2m+3}))=z_{m+1}$, and $m_{2m+1,2m+3}=z_{m}f_{2m+3}-z_{m+2}f_{2m+1}$. 

\textbf{Claim 2.} 
\begin{gather*}
z_{m}f_{2m+3}-z_{m+2}f_{2m+1} = \\
=\sum_{i+j=2m+3,\;  i \geq 3} (i-1)z_{m+i}f_{j} - \sum_{j-i=2m+3,\;  i \geq 1}(i+1)z_{m-i}f_{j}. 
\end{gather*}
It is easy to see that \[
in(z_{m}f_{2m+3}-z_{m+2}f_{2m+1})=z_{m-1}z_{m+2}^{2} \geq in(z_{m-i}f_{j}), \; in(z_{m+i}f_{j}) \]
and hence this is a standard expression. Now it immediately follows from Buchberger's criterion that the $f_{j}$'s form a  Gr\"{o}bner basis 
for $I_{d}$, and claim 1 follows.

Now to prove claim 2, choose numbers $a_{i}$ so that the coefficient of $t^{3m+3}$ in \[
(\sum_{i=0}^{d}a_{i}z_{i}t^{i})(\sum_{j=0}^{2d}f_{j}t^{j}) \]
is zero. Its coefficient is just $\sum_{i+j=3m+3}a_{i}z_{i}f_{j}= \sum_{i+\nu + \mu = 3m+3}a_{i}z_{i}z_{\nu}z_{\mu}$. Fix $i$, $\nu$, $\mu$. 
Then the coefficient of $z_{i}z_{\nu}z_{\mu}$ is $a_{i}+a_{\nu}+a_{\mu}$. Hence the conditions for the coefficient of $t^{3m+3}$ to be zero are
\begin{gather*}
a_{i}+a_{\nu}+a_{\mu}=0 \\
i+ \nu + \mu =3m+3 \\
0 \leq i, \; \nu, \; \mu \leq d\\
m \leq d-2
\end{gather*}
This system is easy to solve. Set $i=\nu = \mu = m+1$. Then $a_{m+1}=0$. Moreover $a_{m-i}+a_{m+1}+a_{m+i+2}=0$ and so 
$a_{m+i+2}=-a_{m-i}$. Now $a_{m+3}+a_{m}+a_{m}=0$, and hence $a_{m+3}=-2a_{m}$, $a_{m-1}=2a_{m}$. Similarly it is easy to see that 
$a_{m+i}=-(i-1)a_{m}$ and $a_{m-i}=(i+1)a_{m}$. Set $a_{m}=1$. Then \[
\sum_{j-i=2m+3,\;  i \geq 1}(i+1)z_{m-i}f_{j}-\sum_{i+j=2m+3,\;  i \geq 3} (i-1)z_{m+i}f_{j}- z_{m+2}f_{2m+1}+z_{m}f_{2m+3}=0 \]
and claim 2 follows immediately. Hence by~\cite[Theorem 15.3]{Eis94}, 
\begin{gather*}
\dim_{k} \frac{\textstyle k[z_{0}, \ldots , z_{d}]}{\textstyle (f_{0}, \ldots , f_{2d})}= \dim_{k} 
\frac{\textstyle k[z_{0}, \ldots , z_{d}]}{\textstyle (in(f_{0}), \ldots , in(f_{2d}))} = \\
= \dim_{k} \frac{\textstyle k[z_{0}, \ldots , z_{d}]}{\textstyle (z_{0}^{2}, z_{1}^{2}, \ldots , z_{d}^{2}, z_{0}z_{1}, z_{1}z_{2}, z_{2}z_{3}, 
\ldots, z_{d-1}z_{d} )}. 
\end{gather*}
Let $(I_{d}:z_{d})$ be the saturation of $z_{d}$ in $I_{d}$. Then it is easy to see that there is an exact sequence \[
0 \la S_{d} / (I_{d}:z_{d}) \stackrel{z_{d}}{\la} S_{d}/I_{d} \la S_{d}/(I_{d},z_{d}) \la 0. \]
Moreover, $(I_{d},z_{d})=(I_{d-1},z_{d})$, and it is not difficult to check that $(I_{d}:z_{d})=(z_{d},z_{d-1},I_{d-2})$. Hence 
there is an exact sequence 
\begin{equation}
0 \la A_{d-2} \la A_{d} \la A_{d-1} \la 0. 
\end{equation}
Proceed by induction on $d$. By induction 
\begin{gather}
\dim_{k} A_{d-1} = \sum_{i\geq -2} \left( \begin{array}{c}
d-1-i \\
2+i
\end{array} \right) \\
\dim_{k}A_{d-2}= \sum_{i\geq -2} \left( \begin{array}{c}
d-2-i \\
2+i
\end{array} \right) = 
\sum_{i\geq -1} \left( \begin{array}{c}
d-1-i \\
1+i
\end{array} \right).
\end{gather}
Moreover it is not difficult to check that 
\begin{equation}
\left(\begin{array}{c}
d-i\\
2+i
\end{array} \right) = 
\left(\begin{array}{c}
d-i-1\\
2+i
\end{array} \right) + 
\left(\begin{array}{c}
d-i-1\\
1+i
\end{array} \right). 
\end{equation} 
Counting dimensions in (2) and taking into consideration (3), (4) and (5), part (a) of the theorem follows.

If $ch(k)=2$ and $\E \cong \sheaf_{\pone}(-d)$, $d>0$, then the equations $(1)$ are just $z_{i}^{2}=0$, $i=0, \ldots , d$, 
and hence \[
D \cong \frac{\textstyle k[z_{0}, \ldots , z_{d}] }{\textstyle (z_{0}^{2}, \ldots , z_{d}^{2} ) } \cong 
\frac{\textstyle k[z_{0}]}{\textstyle (z_{0}^{2})} \otimes \cdots \otimes \frac{\textstyle k[z_{d}]}{\textstyle (z_{d}^{2})} \]
and \[
length D = 2^{d+1}. \]
This concludes the proof of the theorem.
\end{proof}
\section{Lower bound}
We will need the next lemma.
\begin{lema}\label{exseq1}
For any $n > 0$, there is a short exact sequence \[
0 \la    \sheaf_{\pone}(-n) \stackrel{\lambda}{\la} \oplus_{i=1}^{n+1} \sheaf_{\pone} \la 
\oplus_{i=1}^{n} \sheaf_{\pone}(1) \la 0. \]
\end{lema}
\begin{proof}
Consider the left exact sequence \[
0 \la k[x,y](-n) \stackrel{\lambda}{\la} \oplus_{i=1}^{n+1} k[x,y] \stackrel{\phi}{\la}
 \oplus_{i=1}^{n+} k[x,y](1) \]
where $\lambda$ is defined by \[
\lambda(f)=(x^{n}f,x^{n-1}yf,\dots,xy^{n-1}f,y^{n}f) \]
and $\phi$ by \[
\phi(f_{1},\ldots,f_{n+1})=
(yf_{1}-xf_{2},yf_{2}-xf_{3},\ldots,yf_{n}-xf_{n+1}). \]
The sequence is clearly left exact and exact in degrees $\geq 1$. This induces the sequence of the lemma.
\end{proof}
The next proposition is the key to get a lower bound.
\begin{proposition}\label{pushout}
With \x, \Y, \Z\ as in theorem~\ref{deformtotrivial}, and in addition suppose that $\mathcal{F}$ fits in an exact sequence of $\sheaf_{Z}$-
modules \[
0 \la \mathcal{F} \stackrel{\lambda}{\la} \mathcal{F^{\prime}} \la A \la 0 \]
with $A$ an ample $\sheaf_{Z}$-module. Then there exist an extension $W$ of $Y$ by $\mathcal{F^{\prime}}$ 
that fits in a commutative diagram \[
\begin{CD}
0 @>>> \mathcal{F} @>i>> \sheaf_{\x} @>>> \sheaf_{Y} @>>> 0\\
&&        @V{\lambda}VV         @VVV     @VVV          &        \\
0 @>>> \mathcal{F^{\prime}} @>>> \sheaf_{W} @>>> \sheaf_{Y} @>>> 0
\end{CD}
\]
and 
\[ length_{[Z]} \hilb(W) \leq length_{[Z]} \hilb(\x). \]
\end{proposition}
In particular the sequence of lemma~\ref{exseq1} will give lower bounds when 
$Z \cong \pone$ and $\mathcal{F} \cong \oplus_{i} \sheaf_{\pone}(-d_{i})$.
\begin{proof} 
First construct W. By assumption, $\sheaf_{X} = \sheaf_{Y} \oplus \mathcal{F}$ as a sheaf of abelian groups, 
and let $f \in \mathrm{H}^{2}_{S}(\sheaf_{Y}, \mathcal{F})$ define it's ring structure. 
Define $\sheaf_{W}=: \sheaf_{Y} \oplus \mathcal{F^{\prime}}$ as a 
sheaf of abelian groups, and multiplication by the map $\lambda \circ f $. Then the natural map \[
(1_{Y},\lambda): \sheaf_{X}=\sheaf_{Y} \oplus \mathcal{F} \la \sheaf_{W}=\sheaf_{Y} \oplus \mathcal{F^{\prime}} \]
is a sheaf of rings homomorphism. The existence of the commutative diagramm stated before is clear from the construction of $\sheaf_{W}$. 
Now assume that $\hilb(W)=SpecD$, for an Artin local $k$-algebra $(D,m_{D})$. 
First we want to show the existence of a morphism \[
\hilb(W) \la \hilb(\x). \]
$D$ corresponds a deformation $Z_{D} \subset W \times SpecD$ of \Z\ in $W$. Then there is a 
morphism \[
 Z_{D} \stackrel{q}{\la} \x \times SpecD \]
flat over $SpecD$. Since $q$ is a closed immersion over the closed point of $SpecD$, and flat over $SpecD$, 
it must be a closed immersion. 
Hence $Z_{D}$ is a deformation of \Z\ in \x\ and hence there is a morphism $\hilb(W) \la \hilb(\x)$. 
The proposition will follow if we can show that this map is a closed immersion. For this we need the 
following simple result from commutative algebra.
\begin{lema}
Let $(A,m_{A}) \stackrel{\phi}{\la} (B,m_{B})$ be a local homomorphism of rings such that $A/ m_{A} = B / m_{B}$ 
and the induced map $m_{A} / m_{A}^{2} \la m_{B} / m_{B}^{2}$ is surjective. Then $\phi$ is also surjective.
\end{lema}
To see this, let $T=coker( \phi )$ considering $\phi$ as a map of A-modules. Tensor with $A/ m_{A}$ to get \[
m_{A} / m_{A}^{2} \la B / m_{A}B \la T / m_{A}T \la 0. \]
Since $m_{A} / m_{A}^{2} \la m_{B} / m_{B}^{2}$ is surjective it follows that $m_{B}=m_{A}B +m_{B}^{2}$ 
and hence from Nakayama's lemma, $m_{B}=m_{A}B$. Hence $T=m_{A}T=m_{B}T$ and again by Nakayama, $T=0$.

Hence in order to show that the morphism $\hilb(W) \la \hilb(\x)$ is a closed immersion, it suffices to 
show that \[
H^{0}(\hilb(\x), \sheaf_{\hilb(\x)}) \la H^{0}(\hilb(W),\sheaf_{\hilb(W)}) \]
is surjective, and by the lemma that \[
m_{\x}/ m_{\x}^{2} \la m_{W} / m_{W}^{2} \]
is surjective, where 
$m_{\x}, m_{W}$ are the maximal ideals of $\sheaf_{\hilb(\x)}, \sheaf_{\hilb(W)}$ at $[Z]$ 
respectively. The latter is equivalent to show that the map on tangent spaces 
\[ T_{[Z]}(\hilb(W)) \la T_{[Z]}(\hilb(\x)) \]
is injective. Moreover, $T_{[Z]}(\hilb(W))=\mathrm{Hom}_{Z}(\I_{Z, W} / \I_{Z , W}^{2} , \sheaf_{Z})$, and $
T_{[Z]}(\hilb(\x))=\mathrm{Hom}_{Z}(\I_{Z , \x} / \I_{Z , \x}^{2} ,\sheaf_{Z}) $. Hence we need injectivity of the map \[
\mathrm{Hom}_{Z}(\I_{Z , W} / \I_{Z , W}^{2} , \sheaf_{Z}) \la 
\mathrm{Hom}_{Z}(\I_{Z , \x} / \I_{Z , \x}^{2} ,\sheaf_{Z}). \]
Now by construction $\sheaf_{\x}= \sheaf_{Y} \oplus \mathcal{F}$, and $\sheaf_{W}= 
\sheaf_{Y} \oplus \mathcal{F^{\prime}}$ as $\sheaf_{Z} $-modules, and by the assumptions, there is an exact sequence $\sheaf_{\x}$-modules \[
0 \la \sheaf_{\x}  \la \sheaf_{W} \la A \la 0. \]
This will give the sequence \[
0 \la \I_{Z, \x} \la \I_{Z , W} \la A  \la 0. \]
Tensor this sequence with $\sheaf_{Z}$ over $\sheaf_{\x}$ to get a commutative diagramm \[
\begin{CD}
\I_{Z , \x} / \I_{Z , \x}^{2} @>>> \I_{Z , W}\otimes_{\sheaf_{\x} } \sheaf_{Z} @>>> A @>>> 0\\
@VVV @VVV @VVV &\\
\I_{Z , \x} / \I_{Z , \x}^{2} @>>> \I_{Z , W} / \I_{Z , W}^{2} @>>> N @>>> 0
\end{CD}
\]
All vertical maps are surjective and hence $N$ is ample and hence by applying 
$\mathrm{Hom}_{Z}(\cdot \; , \sheaf_{Z})$ we get that the map on tangent spaces is indeed injective, 
and the proposition follows immediately.
\end{proof}
The following two lemmas are needed for the rest of the paper. 
\begin{lema}\label{exseq2}
For any $d>0$, the short exact sequence of lemma~\ref{exseq1} gives rise to the exact sequence \[
0 \la \sheaf_{\pone}(-2d) \stackrel{\lambda \otimes \lambda}{\la} \sheaf_{\pone}^{{(d+1)^{2}}} \la 
\sheaf_{\pone}^{{d^{2}}} \oplus \sheaf_{\pone}^{2d}(1) \la 0. \]
\end{lema}
\begin{proof}
Let $N=coker(\lambda \otimes \lambda)$. Then $N$ fits in a commutative diagramm \[
\begin{CD}
&&&& 0 && 0 &&\\
&& && @VVV @VVV &\\
0 @>>> \sheaf_{\pone}(-2d) @>>> \sheaf_{\pone}(-d)^{d+1} @>>> \sheaf_{\pone}(1-d)^{d} @>>> 0\\
&& @VVV @VVV @VVV &\\
0 @>>> \sheaf_{\pone}(-2d) @>{\lambda \otimes \lambda}>>  \sheaf_{\pone}^{{(d+1)}^{2}} @>>> N @>>> 0\\
&& && @VVV @VVV &\\
&&&& \sheaf_{\pone}(1)^{d(d+1)} @= \sheaf_{\pone}(1)^{d(d+1)}  &&\\
&& && @VVV @VVV &\\
&&&& 0 && 0 &&
\end{CD}
\]
From the last colummn it is easy to see that $H^{1}(\pone,N^{\ast})=0$. $N$ is locally free of 
rank $d^{2}+2d$, and from the middle row we see that $\deg N=2d$. Hence if $N=\oplus_{i=1}^{d^{2}+2d}
\sheaf_{\pone}(a_{i})$, then $\sum_{i=1}^{d^{2}+2d} a_{i}=2d$, and $a_{i} \leq 1$. Hence \[
N=\sheaf_{\pone}^{{d^{2}}} \oplus \sheaf_{\pone}(1)^{2d} \]
and the lemma follows.
\end{proof}
\begin{lema}
Let $D=k[\underline{x}]/I$ be an Artin local ring, and $I \subset (\underline{x})^{2}$. Then the number of isomorphism 
classes of square zero extensions of $D$ by $k$ is $\dim_{k} I/mI$, $m=(\underline{x})$.
\end{lema}
\begin{proposition}
Let  $d>0$ and $\mathcal{E}$ a semi-negative locally free sheaf on \pone. Let \x\ be obtained 
by a sequence of two infinitesimal extensions of $\sheaf_{\pone}$-algebras 
\begin{gather*}
0 \la \sheaf_{\pone}(-d) \la \sheaf_{\x_{1}} \la \sheaf_{\pone} \la 0 \\ 
 0 \la \mathcal{E} \la \sheaf_{\x} \la \sheaf_{\x_{1}} \la 0 
\end{gather*}
such that both are split as extensions of $\sheaf_{\pone}$-modules. Then \[
length_{[\pone]} \hilb(\x) \geq -\deg \mathcal{E} + rank (\E ) +d +2. \]
\end{proposition}
\begin{proof}
The idea of the proof is by using proposition~\ref{pushout} to reduce to the case of extensions of \pone\ by $
\sheaf_{\pone}$, which we understand completely. 
Let $\mathcal{E} = \oplus_{i=1}^{n} \sheaf_{\pone}(-d_{i})$, with $d_{i}>0$. Construct $Y_{1}$ by forming a pushout diagramm \[
\begin{CD}
0 @>>> \sheaf_{\pone}(-d) @>>> \sheaf_{\x_{1}} @>>> \sheaf_{\pone} @>>> 0\\
&& @VV{\lambda}V @VV{\phi}V @VVV &\\
0 @>>> \sheaf_{\pone}^{d+1} @>>> \sheaf_{Y_{1}} @>>> \sheaf_{\pone} @>>> 0 
\end{CD} \]
where $\lambda$ is as in lemma~\ref{exseq1} and the construction of $Y_{1}$ is given in proposition~\ref{pushout}. 
Then from proposition~\cite[Proposition 4.1]{Tziolas99}, it follows that \[
Y \cong \pone \times Spec k[x_{1},\ldots,x_{d+1}] / (x_{1},\ldots,x_{d+1})^{2}. \]
Do the same for \x\, i.e construct $Z$ so that there is pushout diagramm \[
\begin{CD}
0 @>>> \mathcal{E} @>>> \sheaf_{\x} @>>> \sheaf_{\x_{1}} @>>> 0\\
&& @VVV @VVV @VVV &\\
0 @>>> \sheaf_{\pone}^{a} @>>> \sheaf_{Z} @>>> \sheaf_{\x_{1}} @>>> 0
\end{CD} \]
where $a=\sum_{i=1}^{n}(d_{i}+1)$.

\textbf{claim}: There is an extension $Y$ of $Y_{1}$ by $\sheaf_{\pone}^{a} $, such that there is a 
pullback diagramm \[
\begin{CD}
0 @>>> \sheaf_{\pone}^{a} @>>> \sheaf_{Z} @>>> \sheaf_{\x_{1}} @>>> 0\\
&& @VVV @VVV @VV{\phi}V &\\
0 @>>> \sheaf_{\pone}^{a} @>>> \sheaf_{Y} @>>> \sheaf_{Y_{1}} @>>> 0
\end{CD} \]
and $\sheaf_{Z} = \sheaf_{Y} \times_{\sheaf_{Y_{1}}} \sheaf_{\x_{1}} $.
Then by using proposition~\ref{pushout}, we get that \[
length_{[\pone]} \hilb(\x) \geq length_{[\pone]} \hilb(Z) \geq length_{[\pone]} \hilb(Y). \]
Since $Y$ is an extension of $Y_{1}$ by $\sheaf_{\pone}^{a}$, then again by~\cite[Proposition 4.1]{Tziolas99} \[
\begin{array}{ccc}
Y \cong \pone \times SpecD & \mbox{and} & \dim_{k}D=a+d+2=\sum_{i=1}^{n}(d_{i}+1)+d+1
\end{array} \]
and the proposition follows immediately.

Now the map $\sheaf_{\x_{1}} \la \sheaf_{Y_{1}}$, induces a commutative diagramm 
\[
\begin{CD}
\mathrm{Hom}_{\pone}(\sheaf_{Y_{1}},\sheaf_{\pone}^{a}) @>{\delta^{0}}>> \mathrm{Hom}_{\pone}(
\sheaf_{Y_{1}}^{\otimes 2}, \sheaf_{\pone}^{a}) @>{\delta^{1}}>> \mathrm{Hom}_{\pone}(
\sheaf_{Y_{1}}^{\otimes 3}, \sheaf_{\pone}^{a})\\
@VVV      @VV{\lambda}V     @VVV\\
\mathrm{Hom}_{\pone}(\sheaf_{\x_{1}},\sheaf_{\pone}^{a}) @>{d^{0}}>> \mathrm{Hom}_{\pone}
(\sheaf_{{\x_{1}}}^{\otimes 2}, \sheaf_{\pone}^{a}) @>{d^{1}}>> \mathrm{Hom}_{\pone}
(\sheaf_{\x_{1}}^{\otimes 3}, \sheaf_{\pone}^{a}) 
\end{CD} 
\]
By lemma~\ref{hochschild} we know that \[
\begin{array}{ccc}
Ex(Y_{1}, \sheaf_{\pone}^{a})= H^{2}_{S}(\sheaf_{Y_{1}},\sheaf_{\pone}^{a}) & \mbox{and} & 
Ex(X_{1}, \sheaf_{\pone}^{a})= H^{2}_{S}(\sheaf_{X_{1}},\sheaf_{\pone}^{a}).
\end{array} \]
From the diagramm above we get a vector space map \[
0 \la N=Ker(\sigma) \la Ex(Y_{1}, \sheaf_{\pone}^{a}) \stackrel{\sigma}{\la} Ex(\x_{1}, 
\sheaf_{\pone}^{a}). \]
The claim will follow if we show that $\sigma$ is surjective. This is only a question about the dimensions
 of the vector spaces appearing in the above sequence.
Since $Y \cong Spec k[x_{1}, \ldots ,x_{d+1}] / (x_{1}, \ldots ,x_{d+1})^{2}$, then by lemma 3.5 we 
get \[
\dim_{k}Ex(Y_{1}, \sheaf_{\pone}^{a})=a \cdot \left( \begin{array}{c}
d+2 \\
2
\end{array} \right) = \frac{\textstyle a(d+1)(d+2)}{\textstyle 2}.
\]
Now $\dim_{k} Ex(X_{1}, \sheaf_{\pone}^{a})=a \dim_{k} Ex(X_{1}, \sheaf_{\pone})= a \dim_{k} 
\mathrm{H}^{2}_{S}(\sheaf_{X_{1}}, \sheaf_{\pone})$. Moreover, there is an exact sequence \[
0 \la Ker(d^{0}) \la \mathrm{Hom}_{\pone}(\sheaf_{X_{1}}, \sheaf_{\pone}) \la Im (d^{0}) \la 0. \]
It is not difficult to see that $Ker(d^{0})$ are just derivations and that \[
\begin{array}{ccc}
Ker(d^{0})=\mathrm{Hom}_{\pone}(\Omega_{{X_{1}} / \pone}, \sheaf_{\pone}) & and & \Omega_{{X_{1}} / \pone} \cong 
 \sheaf_{\pone}(-d).
\end{array} \]
Hence since $\sheaf_{X_{1}} \cong \sheaf_{\pone} \oplus \sheaf_{\pone}(-d)$ as an $\sheaf_{\pone}$-module, it follows that \[
\dim_{k} Im(d^{0})=1. \]
Next we need to find the symmetric maps in $Ker(d^{1})$. That is we need maps \[
f: (\sheaf_{\pone} \oplus \sheaf_{\pone}(-d)) \otimes (\sheaf_{\pone} \oplus \sheaf_{\pone}(-d)) \la \sheaf_{\pone} \]
that satisfy the conditions of lemma~\ref{ringstructure}. Let  $\phi_{i} \in 
\mathrm{Hom}_{\pone}(\sheaf_{\pone}(-2d),\sheaf_{\pone})$, $i=1,\ldots , 2d+1$. Then it is easy to see that the required maps, 
$f_{i}$, are defined by \[
f_{0}((x,m),(x^{\prime},m^{\prime}))=xm^{\prime} \]
and \[
f_{i}((x,m),(x^{\prime},m^{\prime}))=\phi_{i}(xm^{\prime}+x^{\prime}m), \;\; i\geq 1. \]
Hence \[
\dim_{k} Ex(X_{1}, \sheaf_{\pone}^{a})=a(2d+2-1)=a(2d+1). \]
So to show that $\sigma$ is surjective all that is needed is to show that 
$\dim_{k} Ker(\sigma)=ad(d-1)/2$. Now it is not difficult to see that 
$Ker(\sigma)$ consists of the classes $[f]$, of symmetric maps $f$, such that $f-\delta^{0}(h) \in Ker(\lambda) $ for some $
h \in \mathrm{Hom}_{\pone}(\sheaf_{Y_{1}}, \sheaf_{\pone}^{a})$. The map $\lambda$ is induced 
by the map \[
(\sheaf_{\pone} \oplus \sheaf_{\pone}(-d)) \otimes (\sheaf_{\pone} \oplus \sheaf_{\pone}(-d)) \la 
(\sheaf_{\pone} \oplus \sheaf_{\pone}^{a}) \otimes (\sheaf_{\pone} \oplus \sheaf_{\pone}^{a}) \]
defined by \[
 ((s,f),(s^{\prime},f^{\prime})) \longmapsto ((s,(x^{d}f,x^{d-1}yf,\ldots,y^{d}f)),(s^{\prime},(
x^{d}f^{\prime},x^{d-1}yf^{\prime},\ldots,y^{d}f^{\prime}))) \]
where $x,y$ are the coordinates of \pone. Now from lemma~\ref{exseq2}, \[
\dim_{k} Ker \lambda = ad^{2}. \]
On the other hand, the maps $f_{i,j,k} \in \mathrm{Hom}_{\pone}(\sheaf_{Y_{1}} \otimes \sheaf_{Y_{1}} , 
\sheaf_{\pone}^{a}) \; k=1,\dots,a$, defined by \[
f_{i,j,k}((s,(x_{i})),(s^{\prime},(x_{i}^{\prime})))
=(\delta_{k}^{\nu}(x_{i} x_{j}^{\prime}-x_{j} x_{i}^{\prime}),1 \leq \nu \leq a) \] 
which are $ad(d+1)/2, $are clearly contained in $Ker \lambda $. 
But since they are not symmetric and no  combination of them is symmetric, none of them or 
any combination can be contained in $Ker(\sigma)$. Since $\dim_{k} Ker \lambda = ad^{2}$, there are $
ad(d-1)/2$ basis elements left. Now purely for dimension reasons this  must be the dimension of 
$Ker(\sigma)$. 
\end{proof}
\section{An Example}
In this example we are going to construct a surface $S$ with an elliptic singularity $e$, and a smooth 
rational curve $\overline{\C}$ in $S$ passing through $e$, such that if $\I=\I_{\overline{\C},S}$ then
\begin{enumerate}
\item $\I^{(n)} / \I^{(n+1)} \cong \sheaf_{\overline{\C}}(-1) ,\;\; \forall n \geq 1$, and 
\item$\overline{\C}$ does not move. That is $\dim_{[\overline{\C}]} \hilb(S)=0$.
\end{enumerate}
Let $E$ be an elliptic curve and $P_{0} \in E$ a point. Take this point to be the zero element for the 
group structure of $E$. Let \[
\x=\p(\sheaf_{E} \oplus \sheaf_{E}(-P_{0})) \stackrel{\pi}{\la} E \]
 be a ruled surface over $E$.
The projection of $\sheaf_{E} \oplus \sheaf_{E}(-P_{0})$ to $\sheaf_{E}$ 
and $\sheaf_{E}(-P_{0})$ correspond to the sections $F,E$ respectively, 
such that $E^{2}=-1$ and $F^{2}=1$. To see this note that $\sheaf_{\x}(E)=\sheaf_{\x}(1)$. 
Hence $E \sim h$ where $h=c_{1}(\sheaf_{\x}(1))$.
 Hence \[
E^{2}=h^{2}=E \cdot h=\deg_{E}(\sheaf_{\x}(1) \otimes \sheaf_{E})= \deg_{E} \sheaf_{E}(-P_{0}) = -1. \]
On the other hand $\sheaf_{\x}(F)=\pi^{\ast}(\sheaf_{E}(P_{0}) ) \otimes \sheaf_{\x}(1)$. Hence $
F \sim E + \pi^{\ast}(P_{0})$, and so $F^{2}=E^{2}+2E \cdot \pi^{\ast}(P_{0}) = -1+2=1$.

Now let $P \in E$ be a non torsion point. That is $nP \neq 0 \; \forall n \geq 1$ for the group structure
 on the points of $E$, having $P_{0}$ as zero element. Take $D=\pi^{\ast} (P)$. Then $\sheaf_{\x}(D+E) 
\otimes \sheaf_{E} = \sheaf_{\x}(D) \otimes \sheaf_{\x}(E) \otimes \sheaf_{E}  =\sheaf_{\x}(D) \otimes 
\sheaf_{E}(-P_{0})=\sheaf_{E}(P-P_{0})$. Hence since $P$ is not a torsion point, 
$\sheaf_{\x}(D+E) \otimes \sheaf_{E}$ is also not a  torsion sheaf.\\
Now let $R= F \cap D$. Blow up \x\ at $R$. Let \[
\x^{\prime}=Bl_{R} \x \stackrel{q}{\la} \x \]
be the blow up of $X$ at $R$ and let $B$ be the exceptional divisor. Then $q^{\ast}D=C+B$, and $q^{\ast}F=F^{\prime}+B$ where $C$,
$F^{\prime}$ are  the birational transforms of $D$, $F$ respectively.
Let \[
\x^{\prime} \stackrel{p}{\la} S \]
be the contraction of $E$ to an elliptic singularity $e$~\cite{KoMo98}, and let $\overline{\C} = p(\C)$. 
Let ${x}= \C \cap  B$, and $\I=\I_{\overline{\C},S}$.

\textbf{Claim}: $\I^{(n)} / \I^{(n+1)} \cong \sheaf_{\overline{\C}}(-1) ,\;\; \forall n \geq 1$.\\
Take $f \in \I^{(n)}$ generating $\I^{(n)} / \I^{(n+1)}$ at $e$. Then $f$ defines a principal divisor 
in a neighbourhood of $e$. Denote by $f$ the closure of this Cartier divisor in $S$. Then \[
p^{\ast}f=n \C + mE +D^{\prime}.  \]
Then I claim that $m \geq n+1$. To see this look at the intersection \[
0=p^{\ast}f \cdot E = n \C \cdot E + m E^{2} + D^{\prime} \cdot E=n-m + D^{\prime} \cdot E. \]
If $m=n$, then $D^{\prime} \cdot E =0$, and hence $D^{\prime} \cap E = \emptyset$. But this implies that $
n(\C + E) \cdot E = 0$. Since it has a section then $\C +E \mid _{E}$ is torsion. But this is not true from
 the choice of \C. Hence $m \geq n+1$, and thus it is possible to write \[
p^{\ast}f = nC + (n+1)E + D^{\prime}, \; \; k \geq 0, \; \; D^{\prime} \geq 0.  \]
Now since $(nC+(n+1)E) \cdot E =-1$, $(kE+D^{\prime}) \cdot E = 1$, and $p^{\ast}f$ has a section at E, we see 
that $H^{0}(\sheaf_{\x}(kE+D^{\prime}) \otimes \sheaf_{E}) \neq 0$. Since $Pic \x ^{\prime} = 
q^{\ast} Pic \x \oplus \mathbf{Z} B$, we see that \[
kE+D^{\prime}=\nu E + q^{\ast} \pi ^{\ast} (D_{E}) + \mu  B \; ,  \; D_{E} \in PicE .\]
Now $\deg \sheaf_{\x}(kE+D^{\prime}) \otimes \sheaf_{E} = 1$. Hence $\sheaf_{\x}(kE+D^{\prime}) \otimes 
\sheaf_{E} = \sheaf_{E}(Q) , \; Q \in E$. Hence $D_{E} \sim \nu P_{0} +Q$. Let $D_{n}=:q^{\ast} 
\pi ^{\ast} Q $. Then 
\begin{equation} 
p^{\ast}f = [ nC +(n+1)E+D_{n} ] + \nu D_{0} + \nu E + \mu B 
\setcounter{equation}{0}
\end{equation}
where $D_{0}= \pi ^{\ast} (P_{0})$. But $(nC +(n+1)E+D_{n}) \cdot E = 0$, and $(\nu D_{0} \nu E + \mu B ) 
\mid _{E} = 0$. Hence \[
\sheaf_{\x}(nC +(n+1)E+D_{n}) \otimes \sheaf_{E} \cong \sheaf_{E}.  \]
From $(6.1)$, by pulling back on $E$, we see that $0 \sim nP+(n+1)(-P_{0})+Q+ \nu P_{0} - \nu P_{0}$, and 
hence $Q \sim (n+1)P_{0} - nP$. But this, and using that $q^{\ast}F=F^{\prime}+B \sim E +D_{0}$ gives that 
\begin{equation}
\sheaf_{\x ^{\prime}}(nC+(n+1)E+D_{n}) \cong \sheaf_{\x ^{\prime}}((n+1)F^{\prime}+B).               
\end{equation}
Apply $p_{\ast}$ to get 
\begin{equation}
n \overline{C} + \overline{D}_{n} \sim (n+1)F^{\prime} +B. 
\end{equation}
From $(6.2)$, $p^{\ast}f$ extends to a section of $\sheaf_{\x ^{\prime}}((n+1)F^{\prime}+B)$, and 
hence this will generate $\I^{(n)} / \I^{(n+1)}$ everywhere except $x$. So if $\I^{(n)} / \I^{(n+1)} = 
\sheaf_{\overline{C}}(W)$, then $supp(W)=\{x \}$. From $(6.3)$ we see that \[
\I^{(n)} \I_{{D_{n}}} = \sheaf_{\overline{C}}(-n \overline{C} - \overline{D}_{n}) = \sheaf_{\overline{C}}(
-(n+1)F^{\prime}-B)=\sheaf_{\overline{C}}(-x)=\sheaf_{\overline{C}}(-1). \]
Now there is an injection \[
0 \la \I^{(n)} \I_{{D_{n}}} / \I \I^{(n)} \I_{{D_{n}}}  \la \I^{(n)} / \I^{(n+1)}. \]
Since this is an isomorphism away from $e$, and the support of the corresponding divisors is $\{ x \}$, it 
has to be an isomorphism and hence \[
\I^{(n)} / \I^{(n+1)} \cong \sheaf_{\overline{C}}(-1) \]
and part 1 follows.

To show that $\overline{\C}$ does not move in $S$, we need the following simple result 
\begin{lema}
Let $f : \x \la Y$, be a proper flat morphism of schemes of finite type over a field k, with $Y$ integral. 
Let $y_{0} \in Y$ be a closed point such that the fiber $\x_{{y_{0}}}$ is irreducible and smooth. 
Then there is a neighborhood $U$ of $y_{0}$ in $Y$ such that $\forall y \in U$, the fiber $\x_{y}$ is irreducible and smooth.
\end{lema}
\begin{proof}
The function $\phi (y)= H^{0}(\x_{y} , \sheaf_{\x_{y} })$ is upper semicontinuous. Since $\phi$ is 
surjective and $\phi (y_{0})=1$, there is a neighbourhood $U_{0}$ , of $y_{0}$ so that $\phi (y)=1$, 
$\forall \; y \in U_{0}$, and hence $\x_{y}$ is connected. Smoothness is an open condition, hence since 
$\x_{{y_{0}}}$ is smooth, $f$ is smooth at all $x \in \x_{{y_{0}}}$. Hence there is $x \in V \subset \x$ 
open so that $f$ is smooth in $V$. Since $f$ is flat, it is also open and hence $U=f(V) \cup U_{0}$ is 
open. Then for all $y \in U$, $\x_{y}$ is nonsingular and connected and hence irreducible.
\end{proof}
Now assume that $\overline{\C}$ moves in $S$. Hence there is a family $\mathcal{C} \subset S \times T$ 
flat over $T$, T irreducible, and $\mathcal{C}_{{t_{0}}} = \overline{\C}$, for a point $t_{0} \in T$. Next we are going to show that: 
\begin{enumerate}
\item The general member of the family(which is nonsingular by the lemma) does not intersect $F^{\prime}$.
\item The general member of the family meets B at exactly one point.
\end{enumerate}
Let $\psi : \mathcal{C} \la T$. Consider the restriction $\psi :(F^{\prime} \times T ) \cap \mathcal{C} 
\la T$. $\psi$ is clearly proper. Since $ \overline{\C} \cap F^{\prime} = \emptyset$, it follows that 
$\psi$ is not surjective. Hence $V=T- \psi ((F^{\prime} \times T ) \cap \mathcal{C} )$ is non empty and 
open. Hence $\forall \; t \in V$, $\mathcal{C}_{t} \cap F^{\prime} = \emptyset$, and 1. follows. 
To see 2. consider $\psi (B \times T) \cap \mathcal{C} \la T$. For the same $V$ as before, 
$\psi ^{-1}(V) \la V$ must be a finite morphism. If not then $B$ must be a component of a fiber and 
hence it meets $F^{\prime}$, which is not possible by the choice of $V$. Since $t_{0} \in V$ and $
\overline{\C} \cap B =\{ x \}$, there must be a neighborhood $U$ of $t_{0}$, such that for all 
$t \in U$, $\mathcal{C}_{t} \cap B$ is one point. And this shows 2.

Now let $\mathcal{C}^{\prime} \subset X^{\prime} \times T$ be the birational transform of $\mathcal{C}$. We may 
assume that $T$ is smooth one-dimensional, and hence we get a family of curves in $X^{\prime}$.

Let $C_{t}$ be the general member of $\mathcal{C}$ such that $C_{t} \cap F^{\prime} = \emptyset$ and $C_{t} \cap B$ is a point. 
Since $Pic X^{\prime} \cong q^{\ast} Pic X \oplus \mathbb{Z} B$, there are numbers a,b,c, so that \[
C^{\prime}_{t} \sim aE+bC+cB \]
Clearly, $C^{\prime}_{t} \neq E , \;\; C, \;\; B$. Moreover, by the choice of $C^{\prime}_{t}$, $C^{\prime}_{t} \cdot B =1$ 
and hence $b-c=1$. Also, $C^{\prime}_{t} \cdot F^{\prime}=0$ and hence $c=0$ and $b=1$. Since $C^{\prime}_{t}$ is irreducible, 
$C^{\prime}_{t} \cdot C \geq 0$ and $C^{\prime}_{t} \cdot E \geq 0$. Hence $a-b+c\geq 0$ and $-a+b \geq 0$. 
Hence $a \geq 1 $ and $a \leq 1$ and so $a=1$. Now it immediately follows that \[
C^{\prime}_{t} \sim E+C \]
So $C^{\prime}_{t} \cdot E =0$ which implies that $C^{\prime}_{t} \cap E = \emptyset$. Hence $e \notin C_{t}$, and so $C_{t}$ 
is Cartier. Apply $p_{\ast}$ to see that $\overline{C} \sim C_{t}$ is also Cartier which is clearly impossible. Hence $\overline{C}$ 
does not move. 

Let us now show that it is not even $\mathbb{Q}$-Cartier. Suppose that $nC$ is Cartier for some $n$. 
This means that $\I^{(n)}$ is locally free, and $(nC) \cdot C = - \deg \I^{(n)} / \I^{(n+1)} =1$. But then, $2=(2nC)\cdot C = 
-\deg \I^{(2n)} / \I^{(2n+1)} =1$.

This example shows that a lower bound, other than a constant, for the length of the Hilbert scheme, of a scheme obtained 
by succesive infinitesimal extensions of \pone\ by $\sheaf_{\pone}(-1)$, does not exist.

It would be interesting to know the $length_{[\overline{C}]} \hilb (S)$, for $S$ and $\overline{C}$ as in the previous example. 
Proposition 3.2 predicts $5$ to be the lower bound. If $length_{[\overline{C}]} \hilb (S)=5$, then this is the best possible.

\end{document}